\documentclass[12pt,leqno]{amsart}
\usepackage{amssymb,verbatim,enumerate,ifthen,cite}
\usepackage[mathscr]{eucal}
\usepackage[utf8]{inputenc}
\usepackage[T1]{fontenc}

\def\N{\mathbb{N}}
\def\R{\mathbb{R}}

\def\Z{\mathbb{Z}}
\def\B{\mathscr{B}}
\def\C{\mathscr{C}}
\def\CC{\mathbb{C}}
\def\K{\mathbb{K}}

\def\A{\mathscr{A}}

\def\S{{\mathscr{S}}}

\def\ker{\mathop{\mbox{\rm ker}}\nolimits}
\long\def\comment#1{}

\newtheorem{theorem}{Theorem}[section]
\newtheorem*{theorem*}{Theorem}
\def\Thm#1#2{\ifthenelse{\equal{#1}{*}}{\begin{theorem*}#2\end{theorem*}}
             {\begin{theorem}\label{T#1}#2\end{theorem}}}

\def\thm#1{Theorem~\ref{T#1}}

\newtheorem{proposition}[theorem]{Proposition}
\newtheorem*{proposition*}{Proposition}
\def\Prp#1#2{\ifthenelse{\equal{#1}{*}}{\begin{proposition*}#2\end{proposition*}}
{\begin{proposition}\label{P#1}#2\end{proposition}}}

\newtheorem{corollary}[theorem]{Corollary}
\newtheorem*{corollary*}{Corollary}
\def\Cor#1#2{\ifthenelse{\equal{#1}{*}}{\begin{corollary*}#2\end{corollary*}}
             {\begin{corollary}\label{C#1}#2\end{corollary}}}
\def\cor#1{Corollary~\ref{C#1}}

\newtheorem{lemma}[theorem]{Lemma}
\newtheorem*{lemma*}{Lemma}
\def\Lem#1#2{\ifthenelse{\equal{#1}{*}}{\begin{lemma*}#2\end{lemma*}}
             {\begin{lemma}\label{L#1}#2\end{lemma}}}
\def\lem#1{Lemma~\ref{L#1}}

\theoremstyle{definition}
\newtheorem{remark}[theorem]{Remark}
\newtheorem*{remark*}{Remark}
\def\Rem#1#2{\ifthenelse{\equal{#1}{*}}{\begin{remark}\rm #2\end{remark}}
             {\begin{remark}\label{R#1}\rm #2\end{remark}}}

\newtheorem{example}[theorem]{Example}
\newtheorem*{example*}{Example}
\long\def\Exa#1#2{\ifthenelse{\equal{#1}{*}}{\begin{example*}\rm #2\end{example*}}
             {\begin{example}\label{Ex#1}\rm #2\end{example}}}

\def\eq#1{{\rm(\ref{E#1})}}
\def\Eq#1#2{\ifthenelse{\equal{#1}{*}}
  {\begin{equation*}\begin{aligned}#2\end{aligned}\end{equation*}}
  {\begin{equation}\begin{aligned}\label{E#1}#2\end{aligned}\end{equation}}}

\begin{document}
\begin{flushright}
\end{flushright}
\vspace{5mm}

\date{\today}

\title{Estimates of linear expressions through factorization}

\author[A. H. Ali]{Ali Hasan Ali}
\address[A. H. Ali]{Doctoral School of Mathematical and Computational Sciences, University of Debrecen, H-4002 Debrecen, Pf.\ 400, Hungary}
\email{ali.hasan@science.unideb.hu}
\author[Zs. P\'ales]{Zsolt P\'ales}
\address[Zs. P\'ales]{Institute of Mathematics, University of Debrecen, 
H-4002 Debrecen, Pf.\ 400, Hungary}
\email{pales@science.unideb.hu}


\thanks{The research of the second author was supported by the K-134191 NKFIH Grant.}

\begin{abstract}
The aim of this paper is to establish various factorization results and then to derive estimates for linear functionals through the use of a generalized Taylor theorem. Additionally, several error bounds are established including applications to the trapezoidal rule as well as to a Simpson formula-type rule.
\end{abstract}

\maketitle

\section{Introduction} 

Let $X$, $Y$ and $Z$ be normed spaces over the field $\K$, where $\K$ stands either for the field of real numbers $\R$ or for the field of complex numbers $\CC$. Assume that $A:X\to Y$ and $B:X\to Z$ are given linear maps. To describe the properties of $A$ and to connect it to that of $B$, it could be useful if $A$ admits a factorization $A=C\circ B$, where $C:Z\to Y$ is also linear map. (For instance, if $B$ and $C$ are bounded linear operators, then $\|A\|\leq \|B\|\|C\|$.)

An obvious necessary condition for a decomposition $A=C\circ B$ is that $\ker B\subseteq \ker A$ (where $\ker(\cdot)$ denotes the null space of the corresponding operator). On the other hand, this is also sufficient.
Indeed, suppose that $\ker B\subseteq \ker A$. We define $C$ on $B(X)$ first. For each element $z\in B(X)$, there exists $x\in X$ such that $z=B(x)$ and then we define $C(z):=A(x)$. This definition is correct since if $x'\in X$ also satisfies $z=B(x')$, then $x-x'\in \ker B$ and hence, by the assumption, $x-x'\in \ker A$, which implies that $A(x')=A(x)$. It is easy to see that the map $C:B(X)\to Y$ so defined is linear. If $B(X)=Z$, then it is obvious that $C$ is uniquely determined. If $B(X)$ is a proper subspace, then $C$ can be extended to a linear map $C:Z\to Y$ (using Hamel bases) arbitrarily. However, it is a more important problem to obtain a factorization $A=C\circ B$ in terms of a bounded linear map $C$ provided that $A$ and $B$ are bounded. The next result establishes a sufficient condition for this factorizability.

\Thm{ABC}{Let $X$ and $Z$ be Banach spaces and $Y$ be a normed space over $\K$. Assume that $A:X\to Y$ and $B:X\to Z$ are bounded linear maps such that $\ker B\subseteq \ker A$ and $B(X)=Z$. Then there exists a unique bounded linear map $C:Z\to Y$ such that $A=C\circ B$.}

\begin{proof} In view of the argument above, there exists a uniquely determined linear map $C:Z\to Y$ such 
that $A=C\circ B$ holds. We only have to prove that $C$ is bounded using that $A$ and $B$ are bounded linear maps.

By Banach's Open Mapping Theorem, $B$ is an open map. Therefore, there exists $r>0$ such that $\B_Z\subseteq B(r\B_X)$. (Here $\B_X$ and $\B_Z$ denote the closed unit balls of the spaces $X$ and $Z$, respectively.)

Let $z\in Z$ be a nonzero vector. Then $z/\|z\|\in\B_Z$ and hence there exists $u\in r\B_X$ such that $z/\|z\|=B(u)$. 
Therefore, $z=\|z\|B(u)=B(\|z\|u)$, where $\|u\|\leq r$. This implies that $C(z)=A(\|z\|u)$, thus 
\Eq{*}{
\|C(z)\|=\|A(\|z\|u)\|\leq (r\|A\|)\|z\|.
}
This inequality is obviously true also for $z=0$ and yields that $\|C\|\leq r\|A\|$ and proves the boundedness of $C$.
\end{proof}

The main purpose of this paper is to investigate the problem of factorization in the following setting. Let $I$ denote the compact interval $[a,b]$ and, for $n\geq0$, let $\C_\K^n(I)$ denote the space of $n$ times continuously differentiable $\K$-valued functions (equipped with the norm $\|f\|_{\C^n}:=\|f\|_\infty+\dots+\|f^{(n)}\|_\infty$). (The space $\C_\K^0(I)$ will simply be denoted by $\C_\K(I)$.)  Furthermore, let $\mu$ be a nonzero $\K$-valued and bounded Borel measure on $[a,b]$ throughout this paper. The main goal is to obtain various estimates for the linear functional $\A_\mu:\C_\K(I)\to\K$ defined by
\Eq{*}{
  \A_\mu(f):=\int_{[a,b]} f(x)d\mu(x).
}
In order to construct $n\in\N$ and a linear map $B:\C_\K^n(I)\to\C_\K(I)$ such that $\ker B\subseteq \ker \A_\mu$, we search for exponential polynomials in the kernel of $\A_\mu$. For this aim, let us define the function $\S_\mu:\CC\to\CC$ by
\Eq{*}{
  \S_\mu(\lambda):=\int_{[a,b]} e^{\lambda x}d\mu(x).
}
The function $\S_\mu$ will be termed the \emph{spectral function related to the measure $\mu$}. (In fact, the Laplace transformation of the measure $\mu$ is strongly connected to this function.) Clearly, (due to the boundedness of $\mu$), $\S_\mu$ is an entire function.

The set of zeros of $\S_\mu$ will be denoted by $\Lambda_\mu$ and will be called the \emph{spectral set related to $\mu$}. For $\lambda\in\Lambda_\mu$, let $m\in\N$ be the \emph{multiplicity of the root of $\lambda$}, i.e., the largest number such that
\Eq{*}{
  \S_\mu^{(j)}(\lambda)
  =\int_{[a,b]} x^je^{\lambda x}d\mu(x)=0 \qquad(j\in\{0,\dots,m-1\}).
}
Since $\mu\neq0$, it follows that, for some $\ell\in\N$, $\S_\mu^{(\ell)}(\lambda)\neq0$, and hence the multiplicity of an element $\lambda\in\Lambda_\mu$ is well-defined and will be denoted by $m(\S_\mu,\lambda)$. The above equality shows that, for $\lambda\in\Lambda_\mu$, the exponential polynomials
\Eq{*}{
  x\mapsto x^je^{\lambda x}\qquad(j\in\{0,\dots,m(\S_\mu,\lambda)-1\})
}
belong to the kernel of $\A_\mu$. 

For fixed elements $\lambda_1,\dots,\lambda_k\in\Lambda_\mu$, we consider the polynomial $P:\CC\to\CC$ given by 
\Eq{Pc}{
  P(\lambda)=(\lambda-\lambda_1)^{m_1}\cdots(\lambda-\lambda_k)^{m_k}
  =c_n \lambda^{n}+\dots+c_1 \lambda+c_0,
}
where $1\leq m_i\leq m(\S_\mu,\lambda_i)$ for all $i$ and $n:=m_1+\dots+m_k$. Clearly, $c_n=1$. Then we define the linear differential operator $D_c:\C_\K^n(I)\to\C_\K(I)$ by the formula
\Eq{D}{
 D_c(f):=c_n f^{(n)}+\dots+c_1 f'+c_0 f \qquad(f\in\C^n(I)).
}
It follows from the theory of ordinary differential equations (see \cite{Wal98}) that $D_c$ is a bounded and surjective linear operator and the exponential polynomials
\Eq{expol}{
  x\mapsto x^je^{\lambda_i x}\qquad(i\in\{1,\dots,k\},\, j\in\{0,\dots,m_i-1\})
}
form a fundamental system of solutions of the differential equation $D_c(f)=0$. In other words, the above exponential polynomials span the kernel of $D_c$ and hence $\ker D_c\subseteq \ker \A_\mu$. In view of \thm{ABC}, there exists a unique bounded linear map $C:\C_\K(I)\to\K$ such that $\A_\mu|_{\C^n_\K(I)}=C\circ D_c$. In what follows, we explicitly construct $C$ and then we present several applications to obtain sharp upper bounds for the error terms of quadrature rules. For standard references about error bounds for quadrature rules, we refer to the monographs \cite{Atk89} by Atkinson and \cite{FaiBur98} by Faires and Burden and to the recent papers \cite{BarCerDra02} by Barnett \textit{et al.}, \cite{CeuNeu02} by Cruz-Uribe and Neugebauer, \cite{MasKutHus12} by Masjed-Jamei \textit{et al.}, \cite{Tal06} by Talman and \cite{Uje07} by Ujevi\'{c}.

\section{Auxiliary results: Generalizations of the Taylor theorem}

In this section, we recall the terminology and main results of the paper \cite{AliPal22} by the authors.

Let $c:=(c_0,\dots,c_n)\in\K^{n+1}$, let $D_c$ be defined by \eq{D} and let $\omega_c\in\C^n(\R)$ denote the unique solution of the initial value problem
\Eq{IV}{
  D_c(\omega_c)=0,\qquad \omega_c^{(i)}(0)=\delta_{i,n-1} \quad(i\in\{0,\dots,n-1\}).
} 
The function $\omega_c$ will be called the \emph{characteristic solution} of the differential equation $D_c(\omega)=0$. 
Let $P_c$ denote its characteristic polynomial which is given by
\Eq{P}{
 P_c(\lambda):=c_n \lambda^{n}+\dots+c_1 \lambda+c_0 \qquad(\lambda\in\CC).
}

\Lem{P}{{\rm(\hspace{.35mm}\cite[Lemma 2.3]{AliPal22})} Let $n\in\N$, $c=(c_0,\dots,c_n)\in\K^{n+1}$ with $c_n=1$ and let $\lambda_1,\dots,\lambda_k\in\CC$ be pairwise distinct roots of the characteristic polynomial $P_c$ with multiplicities $m_1,\dots,m_k$, respectively. Then
\Eq{*}{
  \omega_c(t)=\sum_{i=1}^k\sum_{j=0}^{m_i-1}\frac{(P_i^{-1})^{(m_i-1-j)}(\lambda_i)}{(m_i-1-j)!}\cdot\frac{t^je^{\lambda_i t}}{j!},
}
where $P_i(\lambda):=\prod_{\ell\in\{1,\dots,k\}\setminus\{i\}}(\lambda-\lambda_\ell)^{m_\ell}$. }

The following generalization of Taylor's Theorem is going to play a key role in this paper.

\Thm{1}{{\rm(\hspace{.35mm}\cite[Theorem 3.2]{AliPal22})} Let $n\in\N$, $c=(c_0,\dots,c_n)\in\K^{n+1}$ with $c_n=1$.
Then, for all $f\in\C_\K^n(I)$ and $x,a\in I$, we have
\Eq{TF}{
  f(x)=\sum_{j=0}^{n-1}\Bigg( f^{(j)}(a)\sum_{i=0}^{n-1-j}c_{i+j+1}\omega_c^{(i)}(x-a)\Bigg)+\int_a^x D_c(f)(t)\cdot\omega_c(x-t)dt.
}}

We now present several particular cases of the above theorem. For the formulation of these results, for $n\in\N$, $k\in\N_0$ with $k<n$ and $\gamma\in\K$, we define the function $\zeta_{n,k,\gamma}:\R\to\K$ by
\Eq{zeta}{
  \zeta_{n,k,\gamma}(t)
  :=\sum_{i=0}^\infty\frac{\gamma^it^{i(n-k)+n}}{(i(n-k)+n)!}.
}

\Lem{zeta}{{\rm(\hspace{.35mm}\cite[Lemma 2.4]{AliPal22})} Let $n\in\N$, $k\in\N_0$ with $k<n$ and $\gamma\in\K$. Then, function $\zeta_{n,k,\gamma}$ is the (unique) solution of the initial value problem
\Eq{zIVP}{
  \zeta^{(n+1)}=\gamma\zeta^{(k+1)},\qquad
  \zeta^{(i)}(0)=\delta_{i,n}\quad (i\in\{0,\dots,n\}).
}
In addition, if $j\in\{0,\dots,k\}$, then
\Eq{zder}{
  \zeta_{n,k,\gamma}^{(j)}=\zeta_{n-j,k-j,\gamma}.
}
If $\gamma\in\R\setminus\{0\}$, then, for all $t\in\R$,
\Eq{gamma}{
  \zeta_{n,k,\gamma}(t)
  =\begin{cases}
   \sqrt[n-k]{\gamma^{-n}}\zeta_{n,k,1}(\sqrt[n-k]{\gamma}t) 
   & \text{ if } \gamma>0,\\[3mm]
   \sqrt[n-k]{(-\gamma)^{-n}}\zeta_{n,k,-1}(\sqrt[n-k]{-\gamma}t) 
   & \text{ if } \gamma<0.
   \end{cases}
}
We have
\Eq{*}{
 \zeta_{n,k,0}(t)=\frac{t^n}{n!}, \qquad
 \zeta_{n,0,-1}(t)=\sum_{i=1}^\infty\frac{(-1)^{i-1}t^{in}}{(in)!}, \qquad 
 \zeta_{n,0,1}(t)=\sum_{i=1}^\infty\frac{t^{in}}{(in)!},
}}

\Thm{2}{{\rm(\hspace{.35mm}\cite[Theorem 3.2]{AliPal22})} Let $n\in\N$, $k\in\N_0$ with $k<n$, $\gamma\in\K$. Then, for all $f\in\C_\K^n(I)$ and $x,a\in I$, we have
\Eq{TF2}{
  f(x)=\sum_{j=0}^{k-1}f^{(j)}(a)\frac{(x-a)^j}{j!}
  +\sum_{j=k}^{n-1}f^{(j)}(a) \zeta_{n,k,\gamma}^{(n-j)}(x-a)
  +\int_a^x (f^{(n)}(t)-\gamma f^{(k)}(t))\zeta_{n,k,\gamma}'(x-t) dt.
}}

For any continuous function $h:\R\to\R$ let $\rho^+(h)\in[0,+\infty]$ (resp.\ $\rho^-(h)\in[-\infty,0]$) denote the infimum of the positive roots (resp.\  the supremum of the negative roots) of $h$. The subsequent two statements are the mean value theorems related to \thm{1} and \thm{2}, respectively in the real-valued setting.

\Thm{1MV}{{\rm(\hspace{.35mm}\cite[Theorem 4.3]{AliPal22})} Let $n\in\N$, $c=(c_0,\dots,c_n)\in\R^{n+1}$ with $c_n=1$. Then, for all $f\in\C_\R^n(I)$ and $a,x\in I$ with $\rho^-(\omega_c)\leq x-a\leq\rho^+(\omega_c)$, there exists a point $\xi$ between $a$ and $x$ such that 
\Eq{TFMV}{
  f(x)=\sum_{j=0}^{n-1}\Bigg( f^{(j)}(a)\sum_{i=0}^{n-1-j}c_{i+j+1}\omega_c^{(i)}(x-a)\Bigg)+D_c(f)(\xi)\cdot\int_0^{x-a} \omega_c(t)dt.
}}

\Thm{2MV}{{\rm(\hspace{.35mm}\cite[Theorem 4.4]{AliPal22})} Let $n\in\N$, $k\in\N_0$ with $k<n$ and $\gamma\in\R$ and define $\zeta_{n,k,\gamma}$ by \eq{zeta}. Then, for all $f\in\C_\R^n(I)$ and $x,a\in I$ with $\rho^-(\zeta_{n,k,\gamma}')\leq x-a\leq\rho^+(\zeta_{n,k,\gamma}')$, there exists a point $\xi$ between $a$ and $x$ such that
\Eq{2TF2}{
  f(x)&=\sum_{j=0}^{k-1}f^{(j)}(a)\frac{(x-a)^j}{j!}
  +\sum_{j=k}^{n-1}f^{(j)}(a) \zeta_{n,k,\gamma}^{(n-j)}(x-a)
  +(f^{(n)}(\xi)-\gamma f^{(k)}(\xi))\zeta_{n,k,\gamma}(x-a).
}}

\section{Estimating linear functionals}

Our basic factorization theorem is stated as follows.

\Thm{Appl1}{Let $\mu$ be a nonzero bounded $\CC$-valued Borel measure on $[a,b]$, let $\lambda_1,\dots,\lambda_k\in\Lambda_\mu$ and $m_1,\dots,m_k\in\N$ with $m_i\leq m(\S_\mu,\lambda_i)$ for $i\in\{1,\dots,k\}$. Define $c=(c_0,\dots,c_n)\in\CC^{n+1}$  by \eq{Pc} (where $n:=m_1+\dots+m_k$) and the differential operator $D_c:\C^n_\CC([a,b])\to\C_\CC([a,b])$ by \eq{D}. Let $\omega_c\in\C_\CC^n(\R)$ be the characteristic solution of $D_c(\omega)=0$.
Finally, define $g:[a,b]\to\CC$ by
\Eq{g1}{
  g(t):=\int_{[t,b]}\omega_c(x-t) d\mu(x).
}
Then, for all $f\in\C^n_\CC([a,b])$,
\Eq{fact}{
  \A_\mu(f):=\int_{[a,b]}f(x)d\mu(x)
  =\int_{[a,b]} D_c(f)(t)\cdot g(t)dt.
}
In other words, $\A_\mu|_{\C^n_\CC(I)}=C_g\circ D_c$, where $C_g:\C_\CC(I)\to\CC$ is given by
\Eq{*}{
  C_g(h)=\int_{[a,b]} h(t)g(t)dt.
}}

\begin{proof} Due to the assumptions $\lambda_1,\dots,\lambda_k\in\Lambda_\mu$ and $m_i\leq m(\S_\mu,\lambda_i)$ for $i\in\{1,\dots,k\}$, it follows that the exponential polynomials given in \eq{expol} are in the kernel of $\A_\mu$. On the other hand, these functions form a fundamental system of solutions for the differential equation $D_c(\omega)=0$, which implies that $\ker D_c\subseteq \ker \A_\mu$. Since the differential operator has constant coefficients, it follows that, for all $i\geq0$, the function
\Eq{*}{
  \R\ni x\mapsto \omega_c^{(i)}(x-a)
}
is also in the kernel of $D_c$ and hence in the kernel of $\A_\mu$. Therefore, for all $i\geq0$,
\Eq{*}{
  \int_{[a,b]}\omega_c^{(i)}(x-a)d\mu(x)=0.
}
In what follows, let $\chi_S$ denote the characteristic function of any subset $S$ of $[a,b]$. Applying \thm{1}, the above equalities and finally Fubini's Theorem, for all $f\in\C^n_\CC([a,b])$, we obtain
\Eq{*}{
  \A_\mu(f)&=\int_{[a,b]}f(x)d\mu(x)\\
  &=\int_{[a,b]} \bigg(\sum_{j=0}^{n-1}\bigg( f^{(j)}(a)\!\!\sum_{i=0}^{n-1-j}c_{i+j+1}\omega_c^{(i)}(x-a)\bigg)+\int_a^x D_c(f)(t)\cdot\omega_c(x-t)dt\bigg)d\mu(x)\\
  &=\sum_{j=0}^{n-1}f^{(j)}(a)\!\!\sum_{i=0}^{n-1-j}c_{i+j+1}\int_{[a,b]} \omega_c^{(i)}(x-a)d\mu(x)+\int_{[a,b]} \bigg(\int_a^x D_c(f)(t)\cdot\omega_c(x-t)dt\bigg)d\mu(x)\\
  &=\int_{[a,b]} \bigg(\int_a^x D_c(f)(t)\cdot\omega_c(x-t)dt\bigg)d\mu(x)
  =\int_{[a,b]} \int_{[a,b]} \chi_{[a,x]}(t)\cdot D_c(f)(t)\cdot\omega_c(x-t)dt d\mu(x)\\
  &=\int_{[a,b]} \int_{[a,b]} \chi_{[a,x]}(t)\cdot D_c(f)(t)\cdot\omega_c(x-t) d\mu(x)dt
  =\int_{[a,b]} D_c(f)(t)\int_{[t,b]}\omega_c(x-t) d\mu(x)dt\\
  &=\int_{[a,b]} D_c(f)(t)\cdot g(t)dt.
}
This proves \eq{fact}.
\end{proof}

\Thm{Appl2}{Let $\mu$ be a nonzero bounded $\CC$-valued Borel measure on $[a,b]$, let $0\leq k\leq n$ and $\gamma\in\CC$. Assume that
\Eq{ker}{
  \int_{[a,b]}x^id\mu(x)&=0,\qquad (i\in\{0,\dots,k-1\},\\
  \int_{[a,b]}\exp\big(\sqrt[n-k]{\gamma}\exp\big(\tfrac{2j\pi}{n-k}\pmb{i}\big)x\big)d\mu(x)&=0,\qquad
  (j\in\{0,\dots,n-k-1\},
}
where $\sqrt[n-k]{\gamma}$ denotes the root of order $(n-k)$ of $\gamma$ with the smallest nonnegative argument in the interval $[0,2\pi)$.
Define $g:[a,b]\to\CC$ by
\Eq{g2}{
  g(t):=\int_{[t,b]}\zeta'_{n,k,\gamma}(x-t) d\mu(x).
}
Then, for all $f\in\C^n_\CC([a,b])$,
\Eq{fact2}{
  \A_\mu(f)
  =\int_{[a,b]} (f^{(n)}(t)-\gamma f^{(k)}(t))\cdot g(t)dt.
}}

\begin{proof} Consider the differential operator $D_c:\C^n_\CC(I)\to\C_\CC(I)$ given by $D_c(f):=f^{(n)}-\gamma f^{(k)}$. Its characteristic polynomial $P_c$ is given by $P_c(\lambda)=\lambda^n-\gamma\lambda^k=\lambda^k(\lambda^{n-k}-\gamma)$. The roots of this polynomial are $\lambda=0$ with multiplicity $k$ and the roots of order $(n-k)$ of $\gamma$ (with multiplicities equal to $1$). The second group of roots can be written in the form
\Eq{*}{
  \sqrt[n-k]{\gamma}\exp\big(\tfrac{2j\pi}{n-k}\pmb{i}\big)\qquad
  (j\in\{0,\dots,n-k-1\}).
}
Therefore, a fundamental solution system of the differential equation $D_c(f)=0$ can be obtained as
\Eq{*}{
   \big\{x^i\mid i\in\{0,\dots,k-1\}\big\}\cup
   \big\{\exp\big(\sqrt[n-k]{\gamma}\exp\big(\tfrac{2j\pi}{n-k}\pmb{i}\big)x\big)\mid j\in\{0,\dots,n-k-1\}\big\}.
}
Thus, the condition \eq{ker} ensures that $\ker(D_c)\subseteq\ker(\A_\mu)$. Therefore, the polynomial 
\Eq{*}{
  x\mapsto\sum_{j=0}^{k-1}f^{(j)}(a)\frac{(x-a)^j}{j!}
}
(whose degree is at most $k-1$) is in the kernel of $D_c$ and hence, it belongs to $\ker(\A_\mu)$. On the other hand, by equality \eq{zIVP} of \lem{zeta}, $\zeta_{n,k,\gamma}$ solves the differential equation $D_c(f)=0$ and hence, the mapping
\Eq{*}{
x\mapsto\sum_{j=k}^{n-1}f^{(j)}(a) \zeta_{n,k,\gamma}^{(n-j)}(x-a)
}
belongs to the kernel of $D_c$. Consequently, it is also in $\ker(\A_\mu)$. Combining these inclusions, we obtain that
\Eq{*}{
  \int_{[a,b]}\bigg(\sum_{j=0}^{k-1}f^{(j)}(a)\frac{(x-a)^j}{j!} +\sum_{j=k}^{n-1}f^{(j)}(a) \zeta_{n,k,\gamma}^{(n-j)}(x-a)\bigg)d\mu(x)=0.
}
Now applying \thm{2} and following a similar reasoning as in the proof of \thm{Appl1}, we can conclude that \eq{fact2} holds, where the function $g$ is defined by \eq{g2}.
\end{proof}

Before formulating the next result, we recall the definition of the $p$th norm of a Lebesgue measurable function $f:I\to\CC$ and $p\in[1,\infty]$:
\Eq{*}{
  \|f\|_p:=\begin{cases}
           \Big(\int_{[a,b]} |f(t)|^pdt\Big)^{\frac1p}
           & \mbox{if } p\in[1,\infty),\\[2mm]
           \inf\{s\colon |f(t)|\leq s \mbox{ for a.e.\ } t\in I\}
           & \mbox{if } p=\infty.
           \end{cases}
}

\Cor{Appl1}{Under the notations and assumptions of \thm{Appl1}, for all $f\in\C^n_\CC(I)$ and for all $p,q\in[1,\infty]$ with $\frac1p+\frac1q=1$,
\Eq{est}{
  |\A_\mu(f)|\leq \|D_c(f)\|_p\cdot\|g\|_q.
}
If, in addition, $c=(c_0,\dots,c_n)\in\R^{n+1}$, $\mu$ is a real-valued measure, $f\in\C^n_\R(I)$, and both $g$ and $D_c(f)$ are nonnegative (or nonpositive) on $[a,b]$, then
\Eq{ineq}{
  \A_\mu(f)\geq 0.
}
If $g$ and $D_c(f)$ have opposite signs over $[a,b]$, then this inequality reverses.}

\begin{proof}
 Applying the identity \eq{fact} and H\"older's inequality, we get
\Eq{*}{
  |\A_\mu(f)|=\bigg|\int_{[a,b]} D_c(f)(t)\cdot g(t)dt\bigg|
  \leq \big\|D_c(f)\big\|_p\cdot\|g\|_q,
}
which proves \eq{est}. In the real-valued setting, if $D_c(f)\cdot g$ is nonnegative on $[a,b]$, then we can see that \eq{ineq} holds. If $D_c(f)\cdot g$ is nonpositive, then \eq{ineq} is valid with reversed inequality sign.
\end{proof}

\Cor{CTFA}{Let $\mu$ be a nonzero bounded $\CC$-valued Borel measure on $[a,b]$, let $0\leq k\leq n$ and $\gamma\in\CC$. Assume that the equalities in \eq{ker} hold. Define $g:[a,b]\to\CC$ by \eq{g2}.
Then, for all $f\in\C^n_\CC([a,b])$ and $p,q\in[1,\infty]$ with $\frac1p+\frac1q=1$,
\Eq{est+}{
  \bigg|\int_{[a,b]}f(x)d\mu(x)\bigg|\leq \|f^{(n)}-\gamma f^{(k)}\|_p\cdot\|g\|_q.
}
If, in addition, $\gamma\in\R$, $\mu$ is a real-valued measure and both $g$ and $f^{(n)}-\gamma f^{(k)}$ are nonnegative (or nonpositive) on $[a,b]$, then the inequality \eq{ineq} holds. If $g$ and $f^{(n)}-\gamma f^{(k)}$ have opposite signs over $[a,b]$, then this inequality reverses.}

\begin{proof} The statement of this corollary is an immediate consequence of \thm{Appl2} and Hölder's inequality.
\end{proof}

\section{Application to the trapezoidal rule}

The following lemma is very likely well-known, however, for the sake of completeness, we provide its very simple proof.

\Lem{THI}{For all $t\in\R_+$
\Eq{HI}{
  \frac{t}{\sinh(t)}<1<t\coth(t)
}
and, for all $t\in(0,\pi)$,
\Eq{TI}{
  t\cot(t)<1<\frac{t}{\sin(t)}.
}}

\begin{proof}
To prove \eq{HI}, let $t$ be a positive number. By the Taylor series expansions of $\sinh$ and $\cosh$, we get that
\Eq{*}{
  \sinh(t)=\sum_{k=0}^\infty\frac{t^{2k+1}}{(2k+1)!}>\frac{t}{1!}=t,
}
and
\Eq{*}{
  \sinh(t)=\sum_{k=0}^\infty\frac{t^{2k+1}}{(2k+1)!}
  \leq\sum_{k=0}^\infty\frac{t^{2k+1}}{(2k)!}
  =t\sum_{k=0}^\infty\frac{t^{2k}}{(2k)!}=t\cosh(t).
}
These two inequalities directly imply the left and the right hand side inequality in \eq{HI}, respectively.

To verify \eq{TI}, let $t\in(0,\pi)$. Applying the Taylor mean value theorem to the $\sin$ function at $0$, and then using that $\sin$ is positive over $(0,\pi)$, we can conclude that there exists $s\in(0,t)$ such that
\Eq{*}{
  \sin(t)=\sin(0)+\sin'(0)t+\frac{1}{2}\sin''(s)t^2
  =t-\frac{1}{2}\sin(s)t^2<t,
}
which proves the second inequality in \eq{TI}. 

If $t\in[\frac{1}{2}\pi,\pi)$, then $\cot(t)\leq0$, and hence the first inequality in \eq{TI} is obvious. In the case when $t\in(0,\frac{1}{2}\pi)$, we use the Taylor mean value theorem to the $\tan$ function at $0$, which asserts the existence of $s\in(0,t)$ such that
\Eq{*}{
  \tan(t)=\tan(0)+\tan'(0)t+\frac{1}{2}\tan''(s)t^2
  =t+\frac{\sin(s)}{\cos^3(s)}t^2>t,
}
where we used that the functions $\sin$ and $\cos$ are positive over $(0,\frac{1}{2}\pi)$. This inequality implies the first inequality in \eq{TI} also in the case when $t\in(0,\frac{1}{2}\pi)$.
\end{proof}

Given a compact interval $[a,b]$, the classical trapezoidal rule asserts that, for a twice differentiable function $f:[a,b]\to\R$,
\Eq{*}{
   \frac{1}{b-a}\int_a^bf=\frac{f(a)+f(b)}{2}-R_T(f),
}
where the remainder term $R_T(f)$ has various estimates in terms of the norms of the second derivative of $f$ and the length of the interval $[a,b]$. For instance (see \cite[pp. 252--253]{Atk89}),
\Eq{*}{
  |R_T(f)|\leq \frac{(b-a)^2}{12}\|f''\|_\infty.
}
The aim of this section is to establish various further estimates for $R_T(f)$, which is defined by
\Eq{*}{
   R_T(f):=\frac{f(a)+f(b)}{2}-\frac{1}{b-a}\int_a^bf.
}
Observe that, with $\mu:=\frac12(\delta_a+\delta_b)-\nu$ (where $\delta_t$ denotes the Dirac measure concentrated at $t$ and $\nu$ stands for the normalized Lebesgue measure on $[a,b]$), we can obtain that $R_T(f)=\A_\mu(f)$. The corresponding spectral function is given by
\Eq{TS}{
   \S_\mu(\lambda)
   &=\frac{e^{\lambda a}+e^{\lambda b}}{2}
   -\frac{1}{b-a}\int_a^be^{\lambda x}dx \qquad(\lambda\in\CC).
}

\Lem{TS}{Let $\lambda\in\CC$. Then $\lambda$ is a root of the spectral function $\S_\mu$ given by \eq{TS} if and only if $u:=\lambda \frac{b-a}2$ is a fixed point of the tangent hyperbolic function. The multiplicity of $\lambda$ equals $1$ if $\lambda\neq0$ and equals $2$ of $\lambda=0$.}

\begin{proof}
We can see that $\S_\mu(\lambda)=0$ holds trivially if $\lambda=0$, and then $u=0$ is indeed a fixed point of the tangent hyperbolic function. 

An easy computation yields that, for $\lambda\in\CC\setminus\{0\}$,
\Eq{*}{
   \S_\mu(\lambda)
   &=\frac{e^{\lambda b}+e^{\lambda a}}{2}
   -\frac{e^{\lambda b}-e^{\lambda a}}{\lambda(b-a)}\\
   &=e^{\lambda \frac{a+b}2}
   \bigg(\frac{e^{\lambda \frac{b-a}2}+e^{\lambda \frac{a-b}2}}{2}-\frac{e^{\lambda \frac{b-a}2}-e^{\lambda \frac{a-b}2}}{\lambda(b-a)}\bigg)\\
   &=e^{\lambda \frac{a+b}2}
   \bigg(\cosh(\lambda\tfrac{b-a}2)-\frac{2\sinh(\lambda\frac{b-a}2)}{\lambda(b-a)}\bigg).
}
Now, assume that $\lambda\in\CC\setminus\{0\}$ is a solution of the equation $\S_\mu(\lambda)=0$. Then it satisfies
\Eq{*}{
  \cosh(\lambda\tfrac{b-a}2)
  =\frac{2\sinh(\lambda\frac{b-a}2)}{\lambda(b-a)}.
}
In this case, $\cosh(\lambda\frac{b-a}2)$ cannot be zero, because then $\sinh(\lambda\frac{b-a}2)$ would also be zero, which is impossible. Therefore, 
with the notation $u:=\lambda\frac{b-a}{2}$, the equation $\S_\mu(\lambda)=0$ is equivalent to $u=\tanh(u)$, which shows that $u$ must be a fixed point of the tangent hyperbolic function.
To check the multiplicities of the roots of $\S_\mu$, let $u=\lambda \frac{b-a}2$ be a fixed point of the tangent hyperbolic function. Then $\cosh(\lambda\tfrac{b-a}2)\neq0$ and, for $\lambda\neq0$, we get
\Eq{*}{
 \S'_\mu(\lambda)
   &=\tfrac{a+b}2e^{\lambda \frac{a+b}2}
   \bigg(\cosh(\lambda\tfrac{b-a}2)-\frac{2\sinh(\lambda\frac{b-a}2)}{\lambda(b-a)}\bigg)
   +e^{\lambda \frac{a+b}2}
   \bigg(\cosh(\lambda\tfrac{b-a}2)-\frac{2\sinh(\lambda\frac{b-a}2)}{\lambda(b-a)}\bigg)'\\
   &=e^{\lambda \frac{a+b}2}
   \bigg(\tfrac{b-a}2\sinh(\lambda\tfrac{b-a}2)+\frac{2\sinh(\lambda\frac{b-a}2)}{\lambda^2(b-a)}-\frac{\cosh(\lambda\tfrac{b-a}2)}{\lambda}\bigg)\\
   &=e^{\lambda \frac{a+b}2}\cosh(\lambda\tfrac{b-a}2)
   \bigg(\tfrac{b-a}2\tanh(\lambda\tfrac{b-a}2)+\frac{2\tanh(\lambda\frac{b-a}2)}{\lambda^2(b-a)}-\frac{1}\lambda\bigg)\\
   &=e^{\lambda \frac{a+b}2}\cosh(\lambda\tfrac{b-a}2)
   \bigg(\lambda\big(\tfrac{b-a}2\big)^2+\frac{1}{\lambda}-\frac{1}\lambda\bigg)
   =e^{\lambda \frac{a+b}2}\cosh(\lambda\tfrac{b-a}2)
   \lambda\big(\tfrac{b-a}2\big)^2\neq0,
}
which proves that the multiplicity of $\lambda\neq 0$ equals $1$.

In the case when $\lambda=0$, from the defining formula of $\S_\mu$, we get
\Eq{*}{
  S'_\mu(\lambda)=\frac{ae^{\lambda a}+be^{\lambda b}}{2}-\frac{1}{b-a}\int_a^bxe^{\lambda x}dx,
  \qquad S'_\mu(0)=\frac{a+b}{2}-\frac{1}{b-a}\int_a^bxdx=0,
}
and
\Eq{*}{
  S''_\mu(\lambda)=\frac{a^2e^{\lambda a}+b^2e^{\lambda b}}{2}-\frac{1}{b-a}\int_a^bx^2e^{\lambda x}dx.
}
Thus,
\Eq{*}{
S''_\mu(0)=\frac{a^2+b^2}{2}-\frac{1}{b-a}\int_a^bx^2dx=\frac{a^2+b^2}{2}-\frac{a^2+ab+b^2}{3}
=\frac{(b-a)^2}{6}\neq0,
}
which proves that the multiplicity of $\lambda=0$ equals $2$.
\end{proof}

In order to apply our main theorems to the trapezoidal rule, we shall need to describe the fixed points of the tangent hyperbolic function.

\Lem{FP}{A number $u\in\CC$ is a fixed point of the tangent hyperbolic function, i.e., 
 \Eq{th}{
  \tanh(u)=u
}
holds if and only if $u=v\pmb{i}$, where $v\in\R$ is a fixed point of the tangent function. Furthermore, for all $k\in\Z$, the open interval $\big((k-\frac12)\pi,(k+\frac12)\pi\big)$ contains exactly one fixed point of the tangent function.}

\begin{proof} Assume that $u=v\pmb{i}$, where $v\in\R$ and $\tan(v)=v$. Then
\Eq{*}{
  \tanh(u)=\tanh(v\pmb{i})=\frac{\sinh(v\pmb{i})}{\cosh(v\pmb{i})}
  =\frac{\sin(v)}{\cos(v)}\pmb{i}=\tan(v)\pmb{i}=v\pmb{i}=u,
}
which shows that $u$ is a fixed point of the function $\tanh$.

Assume now that \eq{th} holds for $u=w+v\pmb{i}$, where $w,v\in\R$. Then
\Eq{*}{
  \tanh(u)&=\tanh(w+v\pmb{i})
  =\frac{\tanh(w)+\tanh(v\pmb{i})}{1+\tanh(w)\tanh(v\pmb{i})}
  =\frac{\tanh(w)+\tan(v)\pmb{i}}{1+\tanh(w)\tan(v)\pmb{i}}\\
  &=\frac{(\tanh(w)+\tan(v)\pmb{i})(1-\tanh(w)\tan(v)\pmb{i})}{1+\tanh(w)^2\tan(v)^2}.
}
Thus, \eq{th} is equivalent to
\Eq{*}{
  \frac{\tanh(w)(1+\tan(v)^2)}{1+\tanh(w)^2\tan(v)^2}=w,\qquad
  \frac{\tan(v)(1-\tanh(w)^2)}{1+\tanh(w)^2\tan(v)^2}=v.
}
We show that these two equalities imply that $w=0$ (which then easily implies that $\tan(v)=v$). To the contrary, assume that $w\neq0$. Then $w\neq\tanh(w)$, and hence the first equation implies that $\tan(v)\neq 0$ implying that $v\neq0$. Then
\Eq{*}{
  \frac{w}{v}
  =\frac{\tanh(w)(1+\tan(v)^2)}{\tan(v)(1-\tanh(w)^2)}
  =\frac{\sinh(2w)}{\sin(2v)},
}
equivalently
\Eq{*}{
  \frac{\sin(2v)}{2v}=\frac{\sinh(2w)}{2w}.
}
On the other hand, by \lem{THI}, for every nonzero real numbers $x$ and $y$,
the following inequalities hold
\Eq{*}{
  \frac{\sin(x)}{x}<1<\frac{\sinh(y)}{y}.
}
Therefore, these inequalities with $x=2v$ and $y=2w$ yield an obvious contradiction. This contradiction shows that $w=0$ must hold and hence $\tan(v)=v$.

Let $k\in\Z$. The function $h$ defined by $h(x):=\tan(x)-x$ is differentiable on the open interval $I_k:=\big((k-\frac12)\pi,(k+\frac12)\pi\big)$ with a derivative $h'(x)=\tan^2(x)>0$ for all $x\in I_k\setminus\{k\pi\}$. This shows that $h$ is strictly increasing and hence it can have at most one zero in $I_k$. On the other hand, we have that $\lim_{x\to(k\pm\frac{1}{2})\pi}h(x)=\pm\infty$, which shows that $h$ changes sign over $I_k$ and therefore it has a unique zero in $I_k$.  
\end{proof}

The unique fixed point of the tangent function in the open interval $\big((k-\frac12)\pi,(k+\frac12)\pi\big)$ will be denoted by $\tau_k$ in the sequel.

\Thm{TR}{Let $k\in\N$, $0\leq n_1<\dots<n_k$ be integers, let $a,b\in\R$ with $a<b$ and let 
$\lambda_j:=\frac{2}{(b-a)}\tau_{n_j}$ for $j\in\{1,\dots,k\}$.
Define $(c_0,c_1,\dots,c_{2k})\in\R^{2k+1}$ by the equality
\Eq{*}{
  \big(z^2+\lambda_1^2\big)\cdots\big(z^2+\lambda_k^2\big)
  =c_{2k}z^{2k}+\dots+c_1z^1+c_0=:P_c(z)\qquad (z\in\CC).
}
Then, for all $f\in\C_\K^{2k}([a,b])$,
\Eq{*}{
  \frac{f(a)+f(b)}{2}-\frac{1}{b-a}\int_a^bf
  =\int_a^b D_c(f)(t)\cdot g(t)dt,
}
where 
\Eq{gt}{
  g(t):=
  \begin{cases}
   \displaystyle{\sum_{j=1}^k \frac{\sin(\lambda_j(b-t)/2)\sin(\lambda_j(t-a)/2)}{\lambda_jQ_j(\lambda_j)\sin(\lambda_j(b-a)/2)}} & \mbox{if } n_1>0,\\
   \displaystyle{\frac{(b-t)(t-a)}{2Q_1(0)(b-a)}+\sum_{j=2}^k\frac{\sin(\lambda_j(b-t)/2)\sin(\lambda_j(t-a)/2)}{\lambda_jQ_j(\lambda_j)\sin(\lambda_j(b-a)/2)}} & \mbox{if } n_1=0,
  \end{cases}
}
and $Q_j(z):=\prod_{\ell\in\{1,\dots,k\}\setminus\{j\}}(\lambda_\ell^2-z^2)$ for $j\in\{1,\dots,k\}$.}

\begin{proof} Due to the evenness of the characteristic polynomial $P_c$, we have that $c_1=c_3=\dots=c_{2k-1}=0$ and $c_{2k}=1$. The roots of $P_c$ are  $\lambda_j\pmb{i}$ and $(-\lambda_j\pmb{i})$ for $j\in\{1,\dots,k\}$ with multiplicities equal to $1$ except if $n_1=0$, then $\lambda_1=0$ and its multiplicity is equal to $2$. 

First, consider the case when $0<n_1$. Then all roots of the characteristic polynomial have multiplicity $1$. To apply \lem{P}, observe that the polynomial connected to the root $\lambda_j\pmb{i}$ is given as 
\Eq{*}{
\prod_{\ell\in\{1,\dots,k\}\setminus\{j\}}(z-\lambda_\ell\pmb{i})\prod_{\ell\in\{1,\dots,k\}}(z+\lambda_\ell\pmb{i})=(z+\lambda_j\pmb{i})\prod_{\ell\in\{1,\dots,k\}\setminus\{j\}}(z^2+\lambda_\ell^2)=(z+\lambda_j\pmb{i})Q_j(z\pmb{i}).
}
Similarly, the polynomial connected to the root $(-\lambda_j\pmb{i})$ is given as $(z-\lambda_j\pmb{i})Q_j(z\pmb{i})$. Therefore, in view of \lem{P}, we can obtain the characteristic solution $\omega_c$ of $D_c$ is given by 
\Eq{*}{
  \omega_c(t)
  =\sum_{j=1}^k\bigg(\frac{\exp(\lambda_j t\pmb{i})}{2\lambda_j\pmb{i}Q_j(\lambda_j)}-\frac{\exp(-\lambda_j t\pmb{i})}{2\lambda_j\pmb{i}Q_j(\lambda_j)}\bigg)
  =\sum_{j=1}^k\frac{\sin(\lambda_jt)}{\lambda_j Q_j(\lambda_j)}.
}

In the case when $n_1=0$. Then the first root ($\lambda_1=0$) of the characteristic polynomial has multiplicity $2$ and the other roots have multiplicity $1$. Applying \lem{P} in a similar manner, we can see that the characteristic solution $\omega_c$ of $D_c$ is given by
\Eq{*}{
  \omega_c(t)
  =\frac{t}{Q_1(0)}+\sum_{j=2}^k\frac{\sin(\lambda_jt)}{\lambda_j Q_j(\lambda_j)}.
}

By \thm{Appl1}, $g:[a,b]\to\C$ is defined by formula \eq{g1}. Therefore, for $t\in (a,b]$, we have
\Eq{*}{
  g(t)=\frac{1}{2}\omega_c(b-t)-\frac{1}{b-a}\int_t^b \omega_c(x-t)dx.
}
If $n_1>0$, then
\Eq{*}{
  g(t)&=\sum_{j=1}^k\frac{1}{\lambda_j Q_j(\lambda_j)}\bigg(\frac{1}{2}\sin(\lambda_j(b-t))-\int_t^b \frac{\sin(\lambda_j(x-t))}{b-a}dx\bigg)\\
  &=\sum_{j=1}^k\frac{1}{\lambda_j Q_j(\lambda_j)}\bigg(\frac{1}{2}\sin(\lambda_j(b-t))-\frac{1-\cos(\lambda_j(b-t))}{\lambda_j(b-a)}\bigg)\\
  &=\sum_{j=1}^k\frac{2\cos(\lambda_j(b-t))-2+\lambda_j(b-a)\sin(\lambda_j(b-t))}{2\lambda_j^2 Q_j(\lambda_j)(b-a)}.
}
Then, for $j\in\{1,\dots,k\}$, we have that $\lambda_j(b-a)/2=\tau_{n_j}$, which is a fixed point of the tangent function, therefore, 
\Eq{FP}{
  \lambda_j(b-a)=2\tan(\lambda_j(b-a)/2).
}
Hence
\Eq{*}{
  &\frac{2\cos(\lambda_j(b-t))-2+\lambda_j(b-a)\sin(\lambda_j(b-t))}{2\lambda_j^2 Q_j(\lambda_j)(b-a)}\\
  &\qquad=\frac{\cos(\lambda_j(b-t))-1+\tan(\lambda_j(b-a)/2)\sin(\lambda_j(b-t))}{\lambda_j^2Q_j(\lambda_j)(b-a)}.
}
To simplify the numerator of this fraction, we use the identity, 
\Eq{*}{
  \cos(2u)-1+\tan(v)\sin(2u)
  &=-2\sin(u)^2+2\tan(v)\sin(u)\cos(u)\\
  &=\frac{2\sin(u)}{\cos(v)}(\cos(u)\sin(v)-\sin(u)\cos(v))
  =\frac{2\sin(u)\sin(v-u)}{\cos(v)}
}
with $u:=\lambda_j(b-t)/2$ and $v:=\lambda_j(b-a)/2$. Then, it follows that
\Eq{*}{
  \frac{2\cos(\lambda_j(b-t))-2+\lambda_j(b-a)\sin(\lambda_j(b-t))}{2\lambda_j^2 Q_j(\lambda_j)(b-a)}
  =\frac{\sin(\lambda_j(b-t)/2)\sin(\lambda_j(t-a)/2)}{\lambda_jQ_j(\lambda_j)\sin(\lambda_j(b-a)/2)},
}
which proves the equality \eq{gt} when $n_1>0$.

If $n_1=0$, then \eq{gt} can be obtained analogously. 
\end{proof}

In the particular case when $k=1$, the above theorem simplifies to the following result.

\Cor{TR}{Let $n\in\N\cup\{0\}$, let $a,b\in\R$ with $a<b$ and let 
$\lambda_n:=\frac{2\tau_n}{b-a}$. Then, for all $f\in\C_\K^2([a,b])$,
\Eq{*}{
  \frac{f(a)+f(b)}{2}-\frac{1}{b-a}\int_a^bf
  =\int_a^b (f''+\lambda_n^2 f)(t)\cdot g(t)dt,
}
where 
\Eq{g}{
  g(t):=
  \begin{cases}
   \displaystyle{\frac{\sin(\lambda_n(b-t)/2)\sin(\lambda_n(t-a)/2)}{\lambda_n\sin(\lambda_n(b-a)/2)}} & \mbox{if } n>0,\\
   \displaystyle{\frac{(b-t)(t-a)}{2(b-a)}} & \mbox{if } n=0.
  \end{cases}.
}}

Our final statement is a new error estimate for the trapezoidal rule.

\Thm{TR}{Let $n\in\N\cup\{0\}$, let $a,b\in\R$ with $a<b$ and let 
$\lambda_n:=\frac{2\tau_n}{b-a}$. Then, for all $f\in\C_\K^2([a,b])$,
\Eq{EETR}{
  \bigg|\frac{f(a)+f(b)}{2}-\frac{1}{b-a}\int_a^bf\bigg|
  \leq \begin{cases}
       \frac{1}{12}(b-a)^2\cdot\|f''\|_\infty & \mbox{if } n=0, \\[3mm]
       \frac{(n+1)n\pi}{2\tau_n^3}(b-a)^2\cdot\|f''+\lambda_n^2 f\|_\infty & \mbox{if } n>0, \\[3mm]
       \frac{1}{8}(b-a)\cdot\|f''\|_1 & \mbox{if } n=0,\\[3mm]
       \frac{1+|\cos(\tau_n)|}{4\tau_n|\sin(\tau_n)|}(b-a)\cdot\|f''+\lambda_n^2 f\|_1 & \mbox{if } n>0.
       \end{cases}
}}

\begin{proof} According to \cor{TR}, we have that
\Eq{*}{
  \bigg|\frac{f(a)+f(b)}{2}-\frac{1}{b-a}\int_a^bf\bigg|
  =\bigg|\int_a^b (f''+\lambda_n^2 f)(t)\cdot g(t)dt\bigg|
  \leq \begin{cases}
       \|g\|_1\cdot\|f''+\lambda_n^2 f\|_\infty,\\[2mm]
       \|g\|_\infty\cdot\|f''+\lambda_n^2 f\|_1,
       \end{cases}
}
where $g:[a,b]\to\R$ is given by \eq{g}. Therefore, to verify the first inequality in \eq{EETR}, we have to compute $\|g\|_1$. 

If $n=0$, then $\tau_n=0$, $\lambda_n=0$ and, with the substitution $t=a+s(b-a)$, we get
\Eq{*}{
  \|g\|_1=\int_a^b \frac{(b-t)(t-a)}{2(b-a)} dt
  =\frac12(b-a)^2\int_0^1 (1-s)s ds
  =\frac{1}{12}(b-a)^2.
} 

If $n>0$, then, with the substitution $t=a+s(b-a)$, we get
\Eq{*}{
  \|g\|_1&=\int_a^b \bigg|\frac{\sin(\lambda_n(b-t)/2)\sin(\lambda_n(t-a)/2)}{\lambda_n\sin(\lambda_n(b-a)/2)}\bigg|dt\\
  &=\frac{(b-a)^2}{2\tau_n|\sin(\tau_n)|}\int_0^1 \big|\sin((1-s)\tau_n)\sin(s\tau_n)\big|ds.
}
For this computation, we need to find the zeros of $s\mapsto \sin((1-s)\tau_n)\sin(s\tau_n)$ in $[0,1]$:
\Eq{*}{
  s=\frac{k}{\tau_n}\pi,\qquad s=1+\frac{k}{\tau_n}\pi \qquad(k\in\Z).
}
The $n$th positive fixed point $\tau_n$ of the tangent function is in $[n\pi,(n+\tfrac{1}{2})\pi)$. Then $\frac{k}{\tau_n}\pi\in[0,1]$ if $k\in[0,\frac{\tau_n}{\pi}]$, which holds if and only if $0\leq k\leq n$.
Analogously, $1+\frac{k}{\tau_n}\pi\in[0,1]$ if $k\in[-\frac{\tau_n}{\pi},0]$, which holds if and only if $- n\leq k\leq 0$. Therefore, the roots of the map $s\mapsto \sin((1-s)v)\sin(sv)$ in $[0,1]$ in increasing order are
\Eq{*}{
  s_0:=0<s_1:=1-\frac{n\pi}{\tau_n}
  <s_2:=\frac{\pi}{\tau_n}<s_3:=1-\frac{(n-1)\pi}{\tau_n}<\dots<s_{2 n}:=\frac{n\pi}{\tau_n}<s_{2 n+1}:=1.
}
We have that
\Eq{*}{
  \int\sin((1-s)\tau_n)\sin(s\tau_n)ds
  &=\frac12\int(\cos((1-2s)\tau_n)-\cos(\tau_n))ds \\
  &=-\frac{1}{4\tau_n}\big(\sin((1-2s)\tau_n)+2s\tau_n\cos(\tau_n)\big)=:F(s).
}
Therefore,
\Eq{*}{
  \int_0^1\big|\sin((1-s)v)\sin(sv)\big|ds
  &=\sum_{j=1}^{2 n+1}\int_{s_{j-1}}^{s_{j}}\big|\sin((1-s)v)\sin(sv)\big|ds\\
  &=\sum_{j=1}^{2 n+1}|F(s_{j})-F(s_{j-1})|
  =\bigg|\sum_{j=1}^{2 n+1}(-1)^j(F(s_j)-F(s_{j-1}))\bigg|\\
  &=\bigg|2\sum_{j=0}^{n}F(s_{2j})-F(s_0)-2\sum_{j=0}^{n}F(s_{2j+1})+F(s_{2 n+1})\bigg|.
}
Using that
\Eq{*}{
  F(s_0)=F(0)=-\frac{\sin(\tau_n)}{4\tau_n}=-\frac{\cos(\tau_n)}{4},\qquad
  F(s_{2 n+1})=F(1)=\frac{\sin(\tau_n)}{4\tau_n}-\frac{\cos(\tau_n)}{2}=-\frac{\cos(\tau_n)}{4},
}
and
\Eq{*}{
  \sum_{j=0}^{n}F(s_{2j})
  &=\sum_{j=0}^{n}F\big(\tfrac{j\pi}{\tau_n}\big)
  =-\frac{1}{4\tau_n}\sum_{j=0}^{n}\Big(\sin(\tau_n-2j\pi)+2j\pi\cos(\tau_n)\Big)\\
  &=-\frac{n+1}{4\tau_n}\big(\sin(\tau_n)+ n\pi\cos(\tau_n)\big),\\
  \sum_{j=0}^{n}F(s_{2j+1})
   &=\sum_{j=0}^{n}F\big(1-\tfrac{(n-j)\pi}{\tau_n}\big)
   =\sum_{j=0}^{n}F\big(1-\tfrac{j\pi}{\tau_n}\big)\\
  &=-\frac{1}{4\tau_n}\sum_{j=0}^{n}\Big(\sin(2j\pi-\tau_n)+2(\tau_n-j\pi)\cos(\tau_n)\Big)\\
  &=-\frac{n+1}{4\tau_n}\big(-\sin(\tau_n)+(2\tau_n- n\pi)\cos(\tau_n)\big).
}
Therefore,
\Eq{*}{
  2\sum_{j=0}^{n}F(s_{2j})-F(s_0)&-2\sum_{j=0}^{n}F(s_{2j+1})+F(s_{2 n+1})\\
  &=-\frac{n+1}{\tau_n}\sin(\tau_n)
  +\frac{n+1}{\tau_n}(\tau_n- n\pi)\cos(\tau_n)\\
  &=(n+1)\cos(\tau_n)\Big(-1
  +\frac{1}{\tau_n}(\tau_n- n\pi)\Big)
  =-(n+1) n\frac{\cos(\tau_n)}{\tau_n}\pi.
}
Putting together the pieces, we can conclude that
\Eq{*}{
  \|g\|_1&=\frac{(b-a)^2}{2\tau_n|\sin(\tau_n)|}\int_0^1 \big|\sin((1-s)\tau_n)\sin(s\tau_n)\big|ds\\
  &=\frac{(b-a)^2}{2\tau_n|\sin(\tau_n)|}(n+1) n\frac{|\cos(\tau_n)|}{\tau_n}\pi
  =\frac{(n+1)n\pi}{2\tau_n^3}(b-a)^2.
}
Thus, the second inequality in \eq{EETR} is proved.

To verify the third and fourth inequalities in \eq{EETR}, we need to compute $\|g\|_\infty$.

If $n=0$, then $\lambda_n=0$ and 
\Eq{*}{
  \|g\|_\infty=\sup_{t\in[a,b]}\frac{(b-t)(t-a)}{2(b-a)}
  =\frac{b-a}{8}.
}
If $n>1$, then $\lambda_n>0$ and
\Eq{*}{
  \|g\|_\infty=\sup_{t\in[a,b]}\bigg|\frac{\sin(\lambda_n(b-t)/2)\sin(\lambda_n(t-a)/2)}{\lambda_n\sin(\lambda_n(b-a)/2)}\bigg|.
}
The supremum is attained at a point $t$ belonging to the interior of $[a,b]$, therefore,
\Eq{*}{
  0&=g'(t)=\frac{1}{2\sin(\tau_n)}
  \big(-\cos(\lambda_n(b-t)/2)\sin(\lambda_n(t-a)/2)+\sin(\lambda_n(b-t)/2)\cos(\lambda_n(t-a)/2)\big)\\
  &=\frac{\sin(\lambda_n(a+b-2t)/2)}{2\sin(\tau_n)}.
}
It follows from here that there exists $k\in\Z$ such that
\Eq{*}{
  t=\sigma_k:=\frac{a+b}2-\frac{k\pi}{\lambda_n}.
}
Clearly, $\sigma_k\in (a,b)$ if and only if
\Eq{*}{
   \bigg|\frac{k\pi}{\lambda_n}\bigg|<\frac{b-a}{2},
}
i.e., if and only if
\Eq{*}{
  |k|\pi<\tau_n<(n+\frac12)\pi,
}
which is equivalent to the inequality $|k|\leq n$. 
For the value of $g$ at $\sigma_k$, we have
\Eq{*}{
  g(\sigma_k)
  &=\frac{\sin(\lambda_n(b-\sigma_k)/2)\sin(\lambda_n(\sigma_k-a)/2)}{\lambda_n\sin(\lambda_n(b-a)/2)}
  =\frac{\sin(\frac12\tau_n+k\frac{\pi}{2})\sin(\frac12\tau_n-k\frac{\pi}{2})}{\lambda_n\sin(\tau_n)}\\
  &=\frac{\cos(k\pi)-\cos(\tau_n)}{2\lambda_n\sin(\tau_n)}
  =\frac{(-1)^k-\cos(\tau_n)}{2\lambda_n\sin(\tau_n)}.
}
This shows that the sequence $\big(|g(\sigma_k)|\big)_{k=-n}^n$ takes only the following two values:
\Eq{*}{
  \frac{1-\cos(\tau_n)}{2\lambda_n|\sin(\tau_n)|},\qquad
  \frac{1+\cos(\tau_n)}{2\lambda_n|\sin(\tau_n)|}
}
and hence
\Eq{*}{
  \|g\|_\infty
  =\frac{1+|\cos(\tau_n)|}{2\lambda_n|\sin(\tau_n)|}
  =\frac{1+|\cos(\tau_n)|}{4\tau_n|\sin(\tau_n)|}(b-a).
}
Thus, the proof of the theorem is completed.
\end{proof}

\section{An extension of the Simpson formula}

\Thm{ST}{
Let $a,b\in\R$ with $a<b$, $u=w+\pmb{i}v$, where $w\in\R_+$, $v\in(0,\pi)$ and define $\alpha_u,\beta_u$ by
\Eq{*}{
  \alpha_u:=\frac{\overline{u}\sinh(u)-u\sinh(\overline{u})}{2u\overline{u}(\cosh(u)-\cosh(\overline{u}))},\qquad
  \beta_u:=\frac{u\cosh(u)\sinh(\overline{u})-\overline{u}\cosh(\overline{u})\sinh(u)}{u\overline{u}(\cosh(u)-\cosh(\overline{u}))}.
}
Then, for all $f\in\C_\K^4([a,b])$,
\Eq{*}{
  \bigg|&\alpha_u f(a)+\beta_u f\Big(\frac{a+b}{2}\Big)+\alpha_u f(b)-\frac{1}{b-a}\int_a^bf\bigg|\\
  &\leq 
  \begin{cases}
  \dfrac{(b-a)^3\big(v\sinh(w)-w\sin(v)\big)^2} {32(w^2+v^2)^2wv\sinh(w)\sin(v)}\cdot
  \bigg\|f''''+\dfrac{8(v^2-w^2)}{(b-a)^2}f''+\dfrac{16(w^2+v^2)^2}{(b-a)^4}f\bigg\|_1,\\[5mm]
  \dfrac{(b-a)^4(2\alpha_u+\beta_u-1)}{16(w^2+v^2)^2}
  \cdot \bigg\|f''''+\dfrac{8(v^2-w^2)}{(b-a)^2}f''+\dfrac{16(w^2+v^2)^2}{(b-a)^4}f\bigg\|_\infty.
  \end{cases}
}}

\begin{proof}
We have that
\Eq{chd}{
  \cosh(u)&-\cosh(\overline{u})
  =\cosh(w+\pmb{i}v)-\cosh(w-\pmb{i}v)
  =2\sinh(w)\sinh(\pmb{i}v)=2\pmb{i}\sinh(w)\sin(v),
}
which is different from zero due to the assumptions $0<w$ and $0<v<\pi$ and hence $\alpha_u,\beta_u$ are well-defined. Observe that $\alpha_u,\beta_u$ are invariant with respect to conjugation, therefore, they are real numbers,
and one can also see that $\alpha_u$ and $\beta_u$ are the unique solution of the following system of equations:
\Eq{abu}{
  2\alpha_u \cosh(u)+\beta_u=\frac{\sinh(u)}{u},\qquad
  2\alpha_u \cosh(\overline{u})+\beta_u=\frac{\sinh(\overline{u})}{\overline{u}}.
}
Therefore, each element $z$ of the set $\{u,-u,\overline{u},-\overline{u}\}$ is a root of the function $\varphi:\CC\to\CC$ given by
\Eq{*}{
  \varphi_u(z):=2\alpha_u \cosh(z)+\beta_u-\frac{\sinh(z)}{z}.
}
Now, we construct the measure $\mu_u$ on the Borel subsets of $[a,b]$ by
\Eq{*}{
  \mu_u(A):=\alpha_u\delta_a(A)+\beta_u\delta_{\frac{a+b}2}(A)
  +\alpha_u\delta_b(A)-\frac{1}{b-a}\int_{A}1 \qquad(A\subseteq[a,b]),
}
where $\delta_t$ denotes the Dirac measure concentrated at $t\in[a,b]$, and we consider the linear function $\A_{\mu_u}:\C_\K\to\K$ defined by
\Eq{*}{
  \A_{\mu_u}(f):=\int_{[a,b]}fd\mu_u
  =\alpha_u f(a)+\beta_u f\Big(\frac{a+b}{2}\Big)+\alpha_u f(b)-\frac{1}{b-a}\int_a^bf.
}
It is not difficult to see that
\Eq{*}{
  |\A_{\mu_u}(f)|\leq(2|\alpha_u|+|\beta_u|+1)\|f\|_\infty,
}
and the constant $2|\alpha_u|+|\beta_u|+1$ is the sharpest possible one. Furthermore, if $p\in[1,\infty)$ and $2|\alpha_u|+|\beta_u|>0$, then $\A_{\mu_u}$ is not bounded with respect to the norm $\|\cdot\|_p$. 

The spectral function of $\A_{\mu_u}$ is now of the form
\Eq{*}{
  \S_{\mu_u}(\lambda)
  &=\alpha_u e^{\lambda a}+\beta_u e^{\lambda\frac{a+b}{2}}+\alpha_u e^{\lambda b}-\frac{1}{b-a}\int_a^be^{\lambda x}dx\\
  &=\alpha_u e^{\lambda a}+\beta_u e^{\lambda\frac{a+b}{2}}+\alpha_u e^{\lambda b}-\frac{e^{\lambda b}-e^{\lambda a}}{\lambda(b-a)}\\
  &=e^{\lambda\frac{a+b}{2}}\bigg(2\alpha_u\cosh\Big(\lambda\frac{b-a}{2}\Big)+\beta_u-\frac{2}{\lambda(b-a)}\sinh\Big(\lambda\frac{b-a}{2}\Big)\bigg).
}
Therefore, the equation $\S_{\mu_u}(\lambda)=0$ holds if and only if $\varphi_u\big(\frac12\lambda (b-a)\big)=0$ is valid. Consequently, each element of the set $\{\lambda_u,-\lambda_u,\overline{\lambda_u},-\overline{\lambda_u}\}$, where $\lambda_u:=\frac{2u}{b-a}$, is a root of the spectral function $\S_{\mu_u}$. The polynomial whose roots are these numbers is equal to 
\Eq{*}{
P_{u}(z)=\big(z^2-\lambda_u^2\big)\big(z^2-\overline{\lambda_u}^2\big)=z^4-(\lambda_u^2+\overline{\lambda_u}^2)z^2+\lambda_u^2 \overline{\lambda_u}^2.
}
The differential operator whose characteristic polynomial is $P_{u}$ is given by
\Eq{*}{
  D_u(f)=f''''-(\lambda_u^2+\overline{\lambda_u}^2)f''+\lambda_u^2 \overline{\lambda_u}^2f.
}
The characteristic solution of the differential equation $D_u(f)=0$ is now of the form
\Eq{ou}{
  \omega_u(t)
  &=\frac{\sinh(\lambda_u t)}{\lambda_u(\lambda_u^2-\overline{\lambda_u}^2)}+\frac{\sinh(\overline{\lambda_u}t)}{\overline{\lambda_u}(\overline{\lambda_u}^2-\lambda_u^2)}
  = \sum_{k=1}^\infty\frac{\lambda_u^{2k}-\overline{\lambda_u}^{2k}}{\lambda_u^2-\overline{\lambda_u}^2}\frac{t^{2k+1}}{(2k+1)!}.
}
Finally, we compute
\Eq{*}{
  g_u(t)&=\int_{[t,b]}\omega_u(x-t)d\mu_u(x)
  =\int_{[a,b]}\chi_{[t,b]}(x)\cdot\omega_u(x-t)d\mu_u(x)\\
  &=\begin{cases}
    \beta_u\omega_u\big(\frac{a+b}{2}-t\big) +\alpha_u\omega_u(b-t)-\frac{1}{b-a}\int_0^{b-t}\omega_u(x)dx& \mbox{if } a\leq t\leq \frac{a+b}{2}, \\[2mm]
    \alpha_u\omega_u(b-t)-\frac{1}{b-a}\int_0^{b-t}\omega_u(x)dx& \mbox{if } \frac{a+b}{2}<t\leq b, 
   \end{cases}
}
and, if $a\leq t\leq \frac{a+b}{2}$, then, according to \eq{ou},
\Eq{*}{
 g_u(t)&=\beta_u\bigg(\frac{\sinh(\lambda_u \big(\frac{a+b}{2}-t\big))}{\lambda_u(\lambda_u^2-\overline{\lambda_u}^2)}+\frac{\sinh(\overline{\lambda_u}\big(\frac{a+b}{2}-t\big))}{\overline{\lambda_u}(\overline{\lambda_u}^2-\lambda_u^2)}\bigg)+\alpha_u\bigg(\frac{\sinh(\lambda_u (b-t))}{\lambda_u(\lambda_u^2-\overline{\lambda_u}^2)}+\frac{\sinh(\overline{\lambda_u}(b-t))}{\overline{\lambda_u}(\overline{\lambda_u}^2-\lambda_u^2)}\bigg)\\
 &\quad -\frac{1}{b-a}\bigg(\frac{\cosh(\lambda_u (b-t))}{\lambda_u^2(\lambda_u^2-\overline{\lambda_u}^2)}+\frac{\cosh(\overline{\lambda_u}(b-t))}{\overline{\lambda_u}^2(\overline{\lambda_u}^2-\lambda_u^2)}+\frac{1}{\lambda_u^2\overline{\lambda_u}^2}\bigg).
}
If $\frac{a+b}{2}<t\leq b$, then
\Eq{*}{
 g_u(t)&=\alpha_u\bigg(\frac{\sinh(\lambda_u (b-t))}{\lambda_u(\lambda_u^2-\overline{\lambda_u}^2)}+\frac{\sinh(\overline{\lambda_u}(b-t))}{\overline{\lambda_u}(\overline{\lambda_u}^2-\lambda_u^2)}\bigg)\\
 &\quad -\frac{1}{b-a}\bigg(\frac{\cosh(\lambda_u (b-t))}{\lambda_u^2(\lambda_u^2-\overline{\lambda_u}^2)}+\frac{\cosh(\overline{\lambda_u}(b-t))}{\overline{\lambda_u}^2(\overline{\lambda_u}^2-\lambda_u^2)}+\frac{1}{\lambda_u^2\overline{\lambda_u}^2}\bigg).
}
In what follows, we shall prove the following properties of $g_u$.
\begin{enumerate}[(1)]
 \item $g_u$ is symmetric with respect to the midpoint of $[a,b]$, i.e., $g_u(a+b-t)=g_u(t)$ for all $t\in[a,b]$.
 \item $g_u$ is continuous and nonnegative on $[a,b]$, and $g_u(a)=g_u(b)=0$.
 \item $g_u$ is increasing on $\big[a,\frac{a+b}{2}\big]$ and decreasing on $\big[\frac{a+b}{2},b\big]$, consequently $\|g_u\|_\infty=g_u(\frac{a+b}{2})$.
\end{enumerate}
To prove assertion (1), for $\frac{a+b}{2}\leq t\leq b$, we show that the equality $g_u(a+b-t)=g_u(t)$ holds, that is,
\Eq{*}{
  &\beta_u\bigg(\frac{\sinh\big(\lambda_u \big(t-\frac{a+b}{2}\big)\big)}{\lambda_u}-\frac{\sinh\big(\overline{\lambda_u}\big(t-\frac{a+b}{2}\big)\big)}{\overline{\lambda_u}}\bigg)\\
  &\quad+\alpha_u\bigg(\frac{\sinh(\lambda_u (t-a))}{\lambda_u}-\frac{\sinh(\overline{\lambda_u}(t-a))}{\overline{\lambda_u}}\bigg)-\alpha_u\bigg(\frac{\sinh(\lambda_u (b-t))}{\lambda_u}-\frac{\sinh(\overline{\lambda_u}(b-t))}{\overline{\lambda_u}}\bigg)\\
 &=\frac{1}{b-a}\bigg(\frac{\cosh(\lambda_u (t-a))}{\lambda_u^2}-\frac{\cosh(\overline{\lambda_u}(t-a))}{\overline{\lambda_u}^2}\bigg)
 -\frac{1}{b-a}\bigg(\frac{\cosh(\lambda_u (b-t))}{\lambda_u^2}-\frac{\cosh(\overline{\lambda_u}(b-t))}{\overline{\lambda_u}^2}\bigg).
}
Using the standard identities for the hyperbolic functions, the above equality is equivalent to 
\Eq{*}{
  \beta_u\bigg(&\frac{\sinh\big(\lambda_u \big(t-\frac{a+b}{2}\big)\big)}{\lambda_u}-\frac{\sinh\big(\overline{\lambda_u}\big(t-\frac{a+b}{2}\big)\big)}{\overline{\lambda_u}}\bigg)\\
  &\quad+2\alpha_u\bigg(\frac{\sinh\big(\lambda_u \big(t-\frac{a+b}{2}\big)\big)\cosh(u)}{\lambda_u}-\frac{\sinh\big(\overline{\lambda_u}\big(t-\frac{a+b}{2}\big)\big)\cosh(\overline{u})}{\overline{\lambda_u}}\bigg)\\
 &=\frac{2}{b-a}\bigg(\frac{\sinh\big(\lambda_u \big(t-\frac{a+b}{2}\big)\big)\sinh(u)}{\lambda_u^2}-\frac{\sinh\big(\overline{\lambda_u}\big(t-\frac{a+b}{2}\big)\big)\sinh(\overline{u})}{\overline{\lambda_u}^2}\bigg).
}
To see that this is valid, one should multiply the first and the second equalities in \eq{abu} by $\frac{1}{\lambda_u}\sinh\big(\lambda_u \big(t-\frac{a+b}{2}\big)\big)$ and $\frac{1}{\overline{\lambda_u}}\sinh\big(\overline{\lambda_u} \big(t-\frac{a+b}{2}\big)\big)$ side by side, respectively, and then subtract the two equations so obtained.

The continuity of $g_u$ at every $t\in[a,b]\setminus\big\{\frac{a+b}{2}\big\}$ is obvious, and, by the definition of $g_u$, this function is also continuous from the left at $t=\frac{a+b}{2}$. Due to the symmetry with respect to the midpoint of $[a,b]$, it follows that $g_u$ is also continuous from the right at $t=\frac{a+b}{2}$, and hence it is also continuous at $t=\frac{a+b}{2}$. The endpoint properties $g_u(a)=g_u(b)=0$ are trivial. The nonnegativity follows from the monotonicity properties of $g_u$, which will be established in what follows. Obviously, in view of the symmetry with respect to the midpoint of $[a,b]$, it is sufficient to show that $g_u$ is decreasing on the upper half of the interval $[a,b]$.

With the substitution $t:=s\frac{a+b}{2}+(1-s)b$, where $s\in[0,1]$, we get that 
\Eq{*}{
\lambda_u(b-t)=\frac{2u}{b-a}\Big(b-\Big(s\frac{a+b}{2}+(1-s)b\Big)\Big)=us. 
}
Therefore, applying also the identity \eq{chd}, we get
\Eq{*}{
g_u&\Big(s\frac{a+b}{2}+(1-s)b\Big)
=\frac{(b-a)^3}{16u^2\overline{u}^2}\bigg(\frac{\overline{u}\sinh(u)-u\sinh(\overline{u})}{\cosh(u)-\cosh(\overline{u})} \bigg(\frac{\overline{u}\sinh(us)}{u^2-\overline{u}^2}+\frac{u\sinh(\overline{u}s)}{\overline{u}^2-u^2}\bigg)\\
 &\qquad\qquad\qquad\qquad\qquad -\bigg(\frac{\overline{u}^2\cosh(us)}{u^2-\overline{u}^2}+\frac{u^2\cosh(\overline{u}s)}{\overline{u}^2-u^2} +1\bigg)\bigg)\\
&=\frac{-(b-a)^3h_u(s)}{16|u|^4(u^2-\overline{u}^2)(\cosh(u)-\cosh(\overline{u}))}
=\frac{(b-a)^3h_u(s)}{128(w^2+v^2)^2wv\sinh(w)\sin(v)},
}
where
\Eq{*}{
  h_u(s):&=-(\overline{u}\sinh(u)-u\sinh(\overline{u}))(\overline{u}\sinh(us)-u\sinh(\overline{u}s))\\
  &\qquad+(\cosh(\overline{u})-\cosh(u))(u^2\cosh(\overline{u}s)-\overline{u}^2\cosh(us)+\overline{u}^2-u^2)\\
  &=4w^2\cosh((1-s)w)\sin(v)\sin(sv)
  +4v^2\cos((1-s)v)\sinh(sw)\sinh(w)\\
 &\quad +4wv\sin(v)\cos(sv)\sinh((1-s)w)
 +4wv\sinh(w)\cosh(sw)\sin((1-s)v)\\
 &\quad-8wv\sinh(w)\sin(v).
}
We have
\Eq{*}{
  h_u'(s)&=-\overline{u}u(\overline{u}\sinh(u)-u\sinh(\overline{u}))(\cosh(us)-\cosh(\overline{u}s))\\
  &\qquad+\overline{u}u(\cosh(u)-\cosh(\overline{u}))(u\sinh(\overline{u}s)-\overline{u}\sinh(us))\\
  &=4(w^2+v^2)\Big(v\sinh(w)\sinh(sw)\sin((1-s)v)-w\sin(v)\sin(sv)\sinh((1-s)w)\Big).
}
The property that $g_u$ is strictly decreasing on the upper half of the interval $[a,b]$, is equivalent to the assertion that $h_u$ is strictly increasing on $[0,1]$. To verify this, we prove the inequality $h'_u(s)>0$ for $s\in(0,1)$, which can be rewritten as
\Eq{*}{
  v\sinh(w)\sinh(sw)\sin((1-s)v)>w\sin(v)\sin(sv)\sinh((1-s)w).
}
After using elementary computations, the above inequality turns out to be equivalent to the following one:
\Eq{*}{
  v(\cot(sv)-\cot(v))>w(\coth(sw)-\coth(w)).
}
Actually, we are going to show that, for all $w\in\R_+$, $v\in(0,\pi)$ and $s\in(0,1)$,  
\Eq{Cot}{
  v(\cot(sv)-\cot(v))>\frac{1-s}{s}>w(\coth(sw)-\coth(w)).
}
We first verify the left hand side inequality in \eq{Cot}.
Using \lem{THI}, for $v\in(0,\pi)$, we have that
\Eq{*}{
  \Big(\cot(v)-\frac{v}{\sin^2(v)}\Big)'
  =-\frac{2}{\sin^2(v)}+\frac{2v\cos(v)}{\sin^3(v)}
  =\frac{2}{\sin^2(v)}(v\cot(v)-1)<0,
}
which shows that the map $v\mapsto \cot(v)-\frac{v}{\sin^2(v)}$ is strictly decreasing on $(0,\pi)$. Therefore, for all $v\in(0,\pi)$ and $s\in(0,1)$, we get
\Eq{*}{
  \cot(sv)-\frac{sv}{\sin^2(sv)}>\cot(v)-\frac{v}{\sin^2(v)}.
}
Consequently, for each fixed $s\in(0,1)$, the derivative of the map $v\mapsto v(\cot(sv)-\cot(v))$ is strictly positive on $(0,\pi)$ and hence it is strictly increasing on $(0,\pi)$. Therefore, for all $v\in(0,\pi)$,
\Eq{*}{
  v(\cot(sv)-\cot(v))>\lim_{v\to0+} v(\cot(sv)-\cot(v))=\frac{1-s}s.
}
Secondly, we verify the right hand side inequality in \eq{Cot}.
Using \lem{THI}, for $w\in\R_+$, we have that
\Eq{*}{
  \Big(\coth(w)-\frac{w}{\sinh^2(w)}\Big)'
  =-\frac{2}{\sinh^2(w)}+\frac{2w\cosh(w)}{\sinh^3(w)}
  =\frac{2}{\sinh^2(w)}(w\coth(w)-1)>0,
}
which shows that the map $w\mapsto \coth(w)-\frac{w}{\sinh^2(w)}$ is strictly increasing on $\R_+$. Therefore, for all $v\in\R_+$ and $s\in(0,1)$, we get
\Eq{*}{
  \coth(sw)-\frac{sw}{\sinh^2(sw)}<\coth(w)-\frac{w}{\sinh^2(w)}.
}
Consequently, for each fixed $s\in(0,1)$, the derivative of the map $w\mapsto w(\coth(sw)-\coth(w))$ is strictly negative on $\R_+$ and hence it is strictly decreasing on $\R_+$. Therefore, for all $w\in\R_+$, 
\Eq{*}{
  w(\coth(sw)-\coth(w))<\lim_{w\to0+} w(\coth(sw)-\coth(w))=\frac{1-s}s,
}
which shows that the right hand side inequality in \eq{Cot} is also valid. The inequality \eq{Cot} implies that $h_u'>0$ on $[0,1]$, whence the strict increasingness of $h_u$ over $[0,1]$ follows.

Then, with $s=1$, we obtain that $h_u(1)=4\big(v\sinh(w)-w\sin(v)\big)^2$, therefore,
\Eq{*}{
  \|g_u\|_\infty=g_u\Big(\frac{a+b}{2}\Big)
 =\frac{(b-a)^3\big(v\sinh(w)-w\sin(v)\big)^2} {32(w^2+v^2)^2wv\sinh(w)\sin(v)}.
}
In order to determine $\|g_u\|_1$, we need the following computation. 
\Eq{*}{
  \int_0^1 h_u(s)ds
  &=\int_0^1 \bigg(-(\overline{u}\sinh(u)-u\sinh(\overline{u}))(\overline{u}\sinh(us)-u\sinh(\overline{u}s))\\
  &\qquad+(\cosh(\overline{u})-\cosh(u))(u^2\cosh(\overline{u}s)-\overline{u}^2\cosh(us)+\overline{u}^2-u^2)\bigg)ds\\
  &=\bigg[-(\overline{u}\sinh(u)-u\sinh(\overline{u}))\Big(\frac{\overline{u}}{u}\cosh(us)-\frac{u}{\overline{u}}\cosh(\overline{u}s)\Big)\\
  &\qquad+(\cosh(\overline{u})-\cosh(u))\Big(\frac{u^2}{\overline{u}}\sinh(\overline{u}s)-\frac{\overline{u}}{u}^2\sinh(us)+(\overline{u}^2-u^2)s\Big)\bigg]_{s=0}^{s=1}\\
  &=\frac{\left(u^2-\bar{u}^2\right)(u \sinh (\bar{u})-\bar{u} \sinh (u)+\bar{u}(\sinh (u)-u) \cosh (\bar{u})+u \cosh (u)(\bar{u}-\sinh (\bar{u})))}{u \bar{u}}\\
  &=\frac{8 v w\big((w \sin (v)+v\sinh(w))(\cosh (w)-\cos (v))-(v^2+w^2)\sinh (w)\sin (v)\big)}{v^2+w^2}.
}
We have
\Eq{*}{
  (v^2+w^2)(2\alpha_u+\beta_u)\sinh(w)\sin(v)
  =(v\sinh(w)+w\sin(v))(\cosh(w)-\cos(v)),
}
Therefore,
\Eq{*}{
  \int_0^1 h_u(s)ds
  &=8 v w(2\alpha_u+\beta_u-1)\sinh(w)\sin(v)
}
and hence
\Eq{*}{
  \int_a^b g_u(x)dx&=2\int_{\frac{a+b}{2}}^bg_u(x)dx
  =(b-a)\int_0^1 g_u\Big(s\frac{a+b}{2}+(1-s)b\Big) ds\\
  &=\frac{(b-a)^4}{128(w^2+v^2)^2wv\sinh(w)\sin(v)}\int_0^1 h_u(s)ds
  =\frac{(b-a)^4(2\alpha_u+\beta_u-1)}{16(w^2+v^2)^2}.
}
Using now \cor{Appl1} with $p=1$ and $p=\infty$, the statement of the theorem follows.
\end{proof}

In the following result, we deduce the Simpson formula with two error terms by taking the limit $u\to0$ in \thm{ST}.

\Cor{ST}{
Let $a,b\in\R$. Then, for all $f\in\C_\K^4([a,b])$,
\Eq{*}{
  \bigg|&\frac{1}{6} f(a)+\frac{2}{3} f\Big(\frac{a+b}{2}\Big)+\frac{1}{6} f(b)-\frac{1}{b-a}\int_a^bf\bigg|
  \leq 
  \begin{cases}
  \dfrac{(b-a)^3}{1152}\cdot
  \big\|f''''\big\|_1,\\[5mm]
  \dfrac{(b-a)^4}{2880}
  \cdot \big\|f''''\big\|_\infty.
  \end{cases}
}}

\begin{proof} In order to derive the statement as the limiting case of \thm{ST} as $u\to0$, we have to compute the limits in the inequality stated by this theorem. An easy computation shows that
\Eq{alp}{
  \alpha_u
  &=\frac{\overline{u}\sum_{n=0}^\infty\frac{u^{2n+1}}{(2n+1)!}-u\sum_{n=0}^\infty\frac{\overline{u}^{2n+1}}{(2n+1)!}}
  {2u\overline{u}\big(\sum_{n=0}^\infty\frac{u^{2n}}{(2n)!}-\sum_{n=0}^\infty\frac{\overline{u}^{2n}}{(2n)!}\big)}  
  =\frac{\sum_{n=1}^\infty\frac{u^{2n}-\overline{u}^{2n}}{(2n+1)!}}
  {2\sum_{n=1}^\infty\frac{u^{2n}-\overline{u}^{2n}}{(2n)!}}
  =\frac{\frac16+\sum_{n=2}^\infty\frac{u^{2(n-1)}+\dots+\overline{u}^{2(n-1)}}{(2n+1)!}}
  {2\big(\frac12+\sum_{n=2}^\infty\frac{u^{2(n-1)}+\dots+\overline{u}^{2(n-1)}}{(2n)!}\big)}.
}
Therefore,
\Eq{*}{
  \lim_{u\to0}\alpha_u=\frac16.
}
Now, using the first equality in \eq{abu}, we get
\Eq{*}{
  \lim_{u\to0}\beta_u
  =\lim_{u\to0}\bigg(\frac{\sinh(u)}{u}-2\alpha_u\cosh(u)\bigg)=1-\frac26\cdot 1=\frac{2}{3}.
}
Similarly,
\Eq{*}{
  \lim_{(w,v)\to(0,0)}&\dfrac{\big(v\sinh(w)-w\sin(v)\big)^2} {(w^2+v^2)^2wv\sinh(w)\sin(v)}
  =\lim_{(w,v)\to(0,0)}\bigg(\dfrac{v\sinh(w)-w\sin(v)} {(w^2+v^2)wv}\bigg)^2\\
 &=\lim_{(w,v)\to(0,0)}\bigg(\dfrac{v\sum_{n=0}^\infty\frac{w^{2n+1}}{(2n+1)!}-w\sum_{n=0}^\infty(-1)^n\frac{v^{2n+1}}{(2n+1)!}} {(w^2+v^2)wv}\bigg)^2
 \\
 &=\lim_{(w,v)\to(0,0)}\bigg(\dfrac{\sum_{n=1}^\infty\frac{w^{2n}-(-1)^nv^{2n}}{(2n+1)!}} {w^2+v^2}\bigg)^2
 =\lim_{(w,v)\to(0,0)}\bigg(\frac16+\dfrac{\sum_{n=2}^\infty\frac{w^{2n}-(-1)^nv^{2n}}{(2n+1)!}} {w^2+v^2}\bigg)^2=\frac1{36}.
}
Using that $\beta_u=\frac{\sinh(u)}{u}-2\alpha_u\cosh(u)=\frac{\sinh(\overline{u})}{\overline{u}}-2\alpha_u\cosh(\overline{u})$ and then applying \eq{alp}, we get
\Eq{*}{
&\lim_{u\to0} \dfrac{2\alpha_u+\beta_u-1}{(w^2+v^2)^2} 
=\lim_{u\to0} \dfrac{2\alpha_u(2-\cosh(u)-\cosh(\overline{u}))+\frac{\sinh(u)}{u}+\frac{\sinh(\overline{u})}{\overline{u}}-2}{2u^2\overline{u}^2}\\
&=\lim_{u\to0} \dfrac{\sum_{n=1}^\infty\frac{u^{2n}+\overline{u}^{2n}}{(2n+1)!}-2\alpha_u\sum_{n=1}^\infty\frac{u^{2n}+\overline{u}^{2n}}{(2n)!}}{2u^2\overline{u}^2}
=\lim_{u\to0} \dfrac{\sum_{n=1}^2\frac{u^{2n}+\overline{u}^{2n}}{(2n+1)!}-2\alpha_u\sum_{n=1}^2\frac{u^{2n}+\overline{u}^{2n}}{(2n)!}}{2u^2\overline{u}^2}\\
&=\lim_{u\to0} \dfrac{20(u^{2}+\overline{u}^{2})+u^{4}+\overline{u}^{4}-2\alpha_u\big(60(u^{2}+\overline{u}^{2})+5(u^{4}+\overline{u}^{4})\big)}{240u^2\overline{u}^2}\\
&=\lim_{u\to0} \Bigg(\dfrac{\big(\frac12+\sum_{n=2}^\infty\frac{u^{2(n-1)}+\dots+\overline{u}^{2(n-1)}}{(2n)!}\big)(20(u^{2}+\overline{u}^{2})+u^{4}+\overline{u}^{4})}{240u^2\overline{u}^2\big(\frac12+\sum_{n=2}^\infty\frac{u^{2(n-1)}+\dots+\overline{u}^{2(n-1)}}{(2n)!}\big)}\\&\hspace{3cm}-\dfrac{\big(\frac16+\sum_{n=2}^\infty\frac{u^{2(n-1)}+\dots+\overline{u}^{2(n-1)}}{(2n+1)!}\big)\big(60(u^{2}+\overline{u}^{2})+5(u^{4}+\overline{u}^{4})\big)}{240u^2\overline{u}^2\big(\frac12+\sum_{n=2}^\infty\frac{u^{2(n-1)}+\dots+\overline{u}^{2(n-1)}}{(2n)!}\big)}\Bigg)\\
&=\lim_{u\to0} \dfrac{\big(\frac12+\frac{u^{2}+\overline{u}^{2}}{4!}\big)(20(u^{2}+\overline{u}^{2})+u^{4}+\overline{u}^{4})-\big(\frac16+\frac{u^{2}+\overline{u}^{2}}{5!}\big)\big(60(u^{2}+\overline{u}^{2})+5(u^{4}+\overline{u}^{4})\big)}{120u^2\overline{u}^2}
=\dfrac{1}{180}.
}
Using the above equalities, upon taking the limit $u\to0$ in \thm{ST}, the result follows.
\end{proof}

\section*{Acknowledgement}

The authors are thankful for the valuable comments and suggestions of the anonymous referee.


\end{document}